\documentclass[dvips]{amsart}
\usepackage{math,epsfig,psfrag}
\usepackage[all]{xy}
\newcommand\BG{\mathsf{BG}}
\newcommand\EG{\mathsf{EG}}
\newcommand\HNN{{\textsf{HNN}}}
\newcommand\A{{\mathcal A}}
\newcommand\M{{\mathcal M}}

\def\pair<#1>{{\ll}#1{\gg}}
\def\9{{\mathbf 1}}
\def\8{{\mathbf 2}}
\def\7{{\mathbf d}}

\begin{document}
\title{Amenability of groups acting on trees}
\date{\today}
\author{Laurent Bartholdi}
\address{Department of Mathematics, Evans Hall, U.C.Berkeley, USA}
\email{laurent@math.berkeley.edu}
\thanks{The author acknowledges support from the ``Swiss National Fund
  for Scientific Research''}
\begin{abstract}
  This note describes the first example of a group that is amenable,
  but cannot be obtained by subgroups, quotients, extensions and
  direct limits from the class of groups locally of subexponential
  growth. It has a balanced presentation
  \[\Delta = \langle b,t|\,[b,t^2]b^{-1},[[[b,t^{-1}],b],b]\rangle.\]
  In the proof, I show that $\Delta$ acts transitively on a
  $3$-regular tree, and that $\Gamma=\langle b,b^{t^{-1}}\rangle$
  stabilizes a vertex and acts by restriction on a binary rooted tree.
  $\Gamma$ is a ``fractal group'', generated by a $3$-state automaton,
  and is a discrete analogue of the monodromy action of iterates of
  $f(z)=z^2-1$ on associated coverings of the Riemann sphere. $\Delta$
  shares many properties with the Thompson group $F$.
  
  I prove briefly some algebraic properties of $\Gamma$, and in
  particular the convergence of quotient Cayley graphs of $\Gamma$
  (aka ``Schreier graphs'') to the Julia set of $f$.
  
  Whenever convenient, the results are expressed in the framework of
  \emph{weakly branch groups}, with extra hypotheses such as
  \emph{contraction}.
\end{abstract}
\maketitle

%%%%%%%%%%%%%%%%%%%%%%%%%%%%%%%%%%%%%%%%%%%%%%%%%%%%%%%%%%%%%%%%
\section{Introduction}
The purpose of this note is twofold: it hints at the connection
between groups acting on trees (\emph{\`a la} Bass-Serre) and groups
acting on rooted trees (\emph{\`a la} Grigorchuk); and it gives a
criterion for amenability and intermediate growth of the latter (and
sometimes the former).

This paper was written in least possible generality that makes the
proofs non-artificial. Many generalizations are possible, and in
particular to the class of ``monomial groups'' defined below.

As a concrete byproduct, the group
\[\Delta = \langle b,t|\,[b,t^2]b^{-1},[[[b,t^{-1}],b],b]\rangle\]
is the first example of a group that is amenable, but cannot be
obtained by subgroups, quotients, extensions and direct limits from
groups locally of subexponential growth (see Theorem~\ref{thm:main}); and it
furthermore has a balanced presentation and acts vertex-transitively
on a $3$-regular tree (see Theorem~\ref{thm:hnn}).

\subsection{Groups of intermediate growth}
Let $G=\langle S\rangle$ be a finitely generated group. Its
\emph{growth function} is $\gamma(n)=\#\{g\in G:g\in S^n\}$. Define a
preorder on growth functions by $\gamma\precsim\delta$ if
$\gamma(n)\le\delta(Cn)$ for some $C\in\N$ and all $n\in\N$, and
denote its symmetric closure by $\sim$. The $\sim$-equivalence class
of $\gamma$ is independent of the choice of $S$. If $\gamma\nsim 2^n$,
then $G$ has \emph{subexponential growth}.  If furthermore
$\gamma\succnsim n^D$ for all $D$, then $G$ has \emph{intermediate
  growth}.  John Milnor asked in 1968~\cite{milnor:5603} whether such
groups existed, and the first example was constructed in the 1980's by
Grigorchuk~\cite{grigorchuk:growth}; see Equation~\eqref{eq:G}.

\subsection{Amenability}
A group is \emph{amenable}~\cite{vneumann:masses} if it admits a
finitely additive invariant measure.  Examples are finite groups and
abelian groups. Amenability is preserved by taking subgroups,
quotients, extensions, and direct limits. The class $\EG$ of
\emph{elementary amenable groups} are those obtained by these
constructions, starting from finite and abelian groups.  Groups of
subexponential growth are also amenable; the elementary amenable
groups of subexponential growth are all of polynomial
growth~\cite{chou:elementary}.

On the other hand, non-abelian free groups are not amenable, and hence
we have a tower \{elementary amenable groups\} $\subseteq$ \{amenable
groups\} $\subseteq$ \{groups with no free subgroup\}.  Mahlon Day
asked in~\cite{day:amen} whether these inclusions are strict. The
first one is, since Grigorchuk's group of intermediate growth is not
elementary amenable. The second one is also
strict~\cite{olshansky:invmean}; for example, the free Burnside group
of exponent $n$ odd at least $665$ is not amenable~\cite{adyan:rw}.

Even in the class of finitely presented groups, both inclusions are
strict: the Grigorchuk group can be embedded in a finitely presented
amenable group\cite{grigorchuk:amenEG}, and Alexander Ol'shanksi\u\i\ 
and Mark Sapir constructed in~\cite{olshansky-s:nonamen}, for all
sufficiently large odd $n$, a non-amenable finitely presented group
satisfying the identity $[X,Y]^n$.

Following Pierre de la Harpe, Rostislav Grigorchuk and Tullio
Ceccherini-Sil\-ber\-stein~\cite[\S~14]{grigorchuk-h-s:paradox}, we denote
by $\BG$ the smallest class of groups containing all groups locally of
subexponential growth\footnote{i.e.\ whose finitely-generated
  subgroups have subexponential growth.}, and closed under taking
subgroups, quotients, extensions and direct limits.  I show in this
note that $\Delta$ is amenable, but does not belong to $\BG$.

\subsection{Groups acting on trees}
Although $\Delta$ is given by a finite presentation, it may also
be defined by an action on the $3$-regular tree. Let $U$ be the binary
rooted tree. Among the many ways the $3$-regular tree $T$ can be
represented, we choose the following two:
\begin{itemize}
\item an infinite horizontal line, called the \emph{axis}, with an
  edge hanging down at each integer coordinate, and a copy of $U$
  attached to that edge's other extremity;
\item a copy of $U$, in which the root vertex has been removed and its
  two adjoining edges have been fused together.
\end{itemize}
The advantage of the first model is that it contains a natural
hyperbolic element, namely the translation one step to the left along
the axis.  For any $n\in\N$, the set of vertices connecting to the
axis at coordinate $\le n$ span a subtree $T_n$ isomorphic to $U$.

\begin{figure}
  \begin{center}
    \psfrag{t}{$t$}
    \psfrag{b}{$b$}
    \psfrag{b1}{$b_0$}
    \psfrag{b2}{$b_0^2$}
    \psfrag{b4}{$b_0^4$} 
    \psfrag{0}{\small $0$}
    \psfrag{-1}{\small $-1$}
    \psfrag{-3}{\small $-3$}
    \psfrag{1}{\small $1$}
    \psfrag{3}{\small $3$}
    \psfrag{5}{\small $5$}
    \psfrag{T2}{\large $T_{-2}$}
    \psfrag{T0}{\large $T_0$}
    \epsfig{file=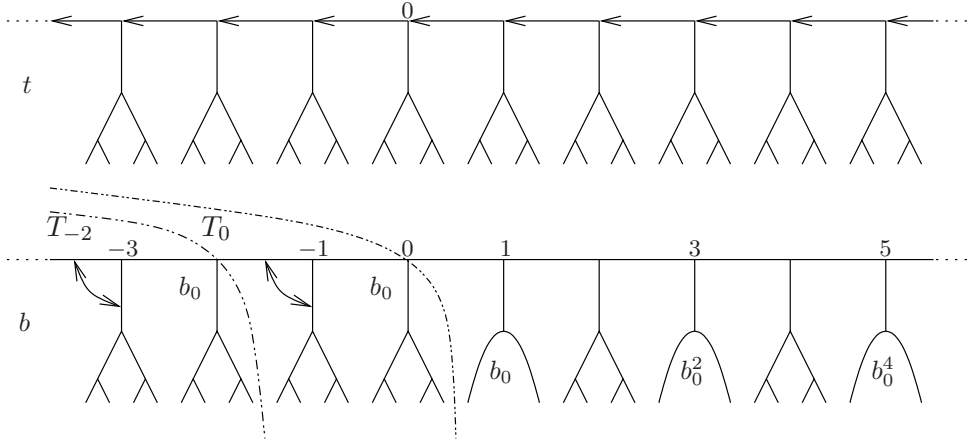}
  \end{center}
  \caption{The action of the generators $b$ and $t$ of $\Delta$ on the
    $3$-regular tree $T$.}
  \label{fig:pinkaction}
\end{figure}

We start by describing the action of $\Delta$ in the first model. Let
$t$ act on $T$ by shifting one step to the left along the axis, and
define the tree isometry $b$ as follows: first, its restriction $b_0$
to the rooted binary tree $T_0$ switches the downward and leftward
branches at $-1-2n$ for every $n\in\N$, starting from $-\infty$ and
moving towards $0$.  Next, identify each of the binary trees below
$n>0$ with $T_0$ in a translation-invariant way. Then $b$ fixes the
half-axis $\N$ and the subtrees below $2n$ for every $n\in\N$, and
acts on the binary tree below $1+2n$ like $b_0^{2^n}$ acts on $T_0$;
see Figure~\ref{fig:pinkaction}.

For convenience, in the sequel, we will always write $b$ for $b_0$ and
$\tilde b$ for $b$.  The action in the second model will be described
in Subsection~\ref{subs:transitive}. Let me just remark that in that
picture $t$ acts as one of the standard generators of the Thompson
group~\cite{cannon-f-p:thompson}.

The stabilizer of $0$ contains $\Gamma=\langle b,b^{t^{-1}}\rangle$,
and $\Delta$ is an ascending \HNN\ extension of $\Gamma$ by $t$.
Indeed writing $a=b^{t^{-1}}$ we have $a^t=b$ and $b^t=a^2$ in
$\Delta$.

The action of $\Gamma$ restricts to a faithful action on $U$, whose
vertices can be naturally labelled by words over $\{\9,\8\}$, with $1$
corresponding to left and horizontal edges and $\8$ corresponding to
right and vertical ones. The action can then be described by
\begin{equation}
  (\9w)^a=\8w^b,\quad(\8w)^a=\9w,\quad(\9w)^b=\9w^a,\quad(\8w)^b=\8w.
  \label{eq:Gamma}
\end{equation}

This is an example of a group \emph{generated by a finite-state
  automaton}. A \emph{transducer} is a tuple $\A=(Q,X,\lambda,\tau)$
with $Q,X$ finite sets called \emph{states} and \emph{letters},
$\lambda:Q\times X\to X$ an \emph{output function} and $\tau:Q\times
X\to Q$ a \emph{transition function}.  A choice of initial state $q\in
Q$ defines an action of $\A_q$ on the tree $X^*$, by
\[()^{\A_q}=(),\quad (xw)^{\A_q}=\lambda(q,x)w^{\A_{\tau(q,x)}}.\]
If each of these transformations is invertible, the \emph{group of
  $\A$} is defined as the group $G(\A)$ generated by $\{\A_q\}_{q\in
  Q}$.

\begin{figure}
  \begin{center}
    \psfrag{0/1}{\footnotesize $0/1$}
    \psfrag{1/0}{\footnotesize $1/0$}
    \psfrag{0/0}{\footnotesize $0/0$}
    \psfrag{1/1}{\footnotesize $1/1$}
    \psfrag{0/1,1/0}{\footnotesize $0/1,1/0$}
    \psfrag{0/0,1/1}{\footnotesize $0/0,1/1$}
    \psfrag{a}{\raisebox{-0.5ex}{$a$}}
    \psfrag{b}{\raisebox{-0.5ex}{$b$}}
    \psfrag{c}{\raisebox{-0.5ex}{$c$}}
    \psfrag{d}{\raisebox{-0.5ex}{$d$}}
    \psfrag{1}{\raisebox{-0.5ex}{$1$}}
    \psfrag{l}[Bc][Bc]{\raisebox{-0.5ex}{$\lambda$}}
    \psfrag{L}[Bc][c]{\raisebox{-0.5ex}{$\lambda^{-1}$}}
    \psfrag{m}[c][Bc]{\raisebox{-0.5ex}{$\mu$}}
    \psfrag{M}[c][c]{\raisebox{-0.5ex}{$\mu^{-1}$}}
    \epsfig{file=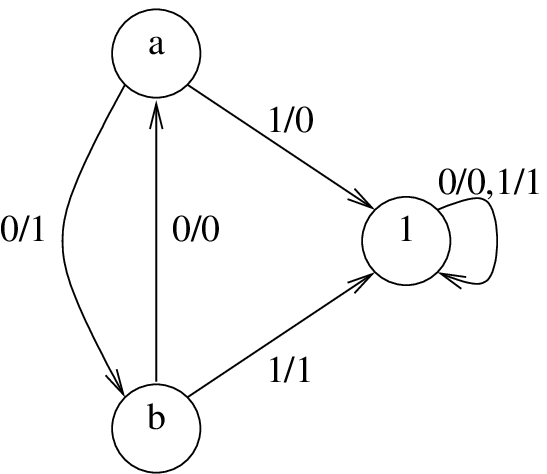}\kern+1em\epsfig{file=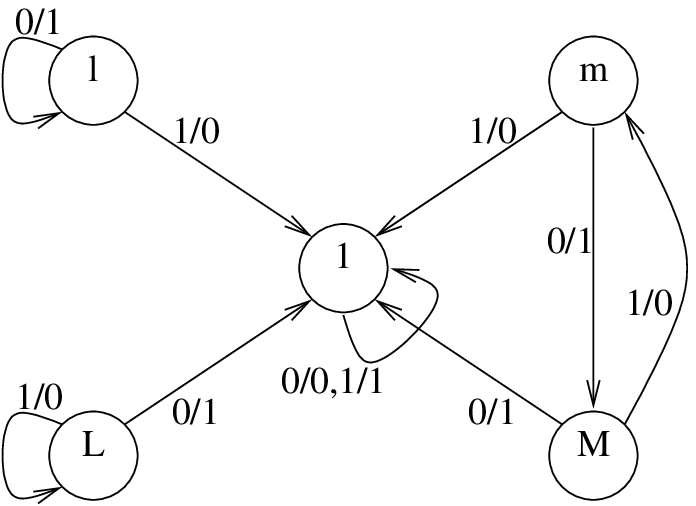}\\
    \epsfig{file=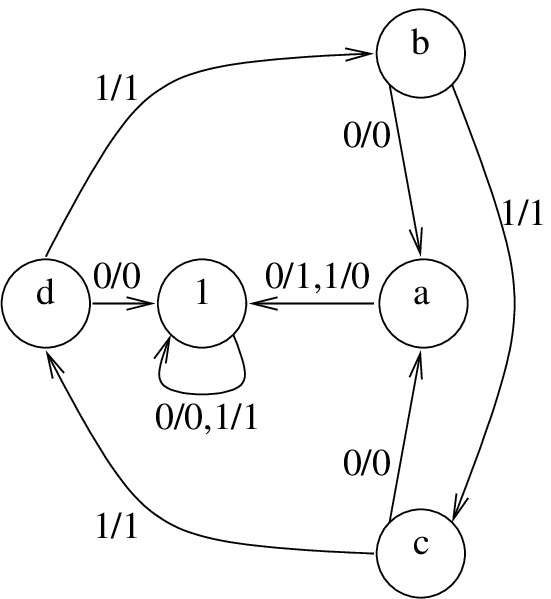}\kern+1em\epsfig{file=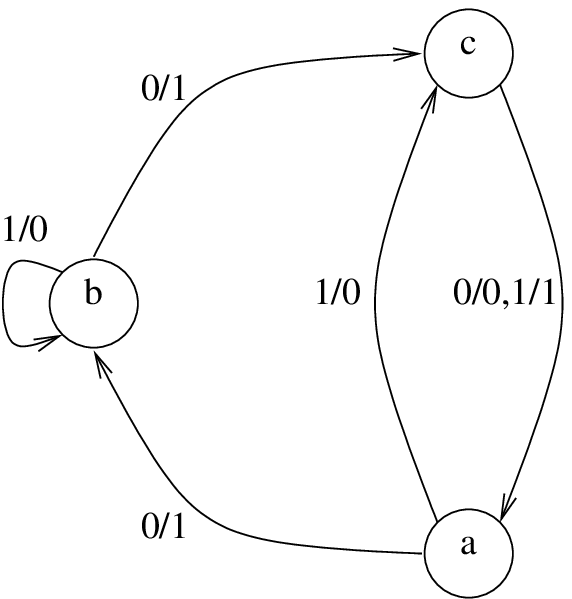}
  \end{center}
  \caption{Automata generating the group $\Gamma$ (top left), the BSV
    group (top right), the Grigorchuk group (bottom left), and the
    Aleshin group (bottom right).}
  \label{fig:autom}
\end{figure}

Automata can be described as graphs, with states as vertices, and an
edge from $q$ to $\tau(q,x)$ labelled $x/\lambda(q,x)$ for all $q\in
Q$ and $x\in X$. Figure~\ref{fig:autom} (top left) gives an automaton
generating $\Gamma$, and Figure~\ref{fig:autom} (bottom left) gives an
automaton generating the Grigorchuk group mentioned above and in
Equation~\eqref{eq:G}.

\subsection{Automata groups}
Automata groups are mainly studied using their \emph{decomposition
  map}: given $g\in G$ acting on the rooted tree $X^*$, its action may
be decomposed in $\#X$ actions on the subtrees connected to the root,
followed by a permutation of the branches at the root. This induces a
group homomorphism $\psi:G\to G\wr\sym X$, written
$\psi(g)=\pair<g_x:x\in X>\pi_g$, into a wreath product\footnote{$\sym
  X$ denotes the symmetric group on $X$; the wreath product is
  $G^X\rtimes\sym X$.}.  The effect of $\psi$ on generators $g\in Q$
can be read quite simply from the automaton: $\pi_g=\lambda(g,-)$ and
$g_x=\tau(g,x)$.

Two favourable situations may occur: first, the definition of an
automaton is dual in that $X,Q$ and $\lambda,\tau$ may be switched
simultaneously. If the dual automaton $\A^*$ generates invertible
transformations of $Q^*$, then the group
\begin{equation}
  \Pi=\langle Q\cup X|\,xq=\tau(q,x)\lambda(q,x)\text{ for all }q\in
  Q,x\in X\rangle
  \label{eq:pi}
\end{equation}
naturally acts on the product of trees $F_X\times F_Q$, and the
original group $G$ can be recovered as the quotient $\langle
Q\rangle/\langle Q\rangle^{X^*}$.

An important example is the automaton $\A$ in Figure~\ref{fig:autom}
(bottom right); it is conjectured that the group it generates is free
on $\A$'s states, though the ``proof'' in~\cite{aleshin:free} appears
to be incomplete.
  
Another favourable situation is the existence of a word metric
$|\cdot|$ on $G$ such that $|g_x|\le\eta|g|+C$ for some $\eta<1$ and
all $g\in G$; this property, called \emph{contraction}, opens the road
to inductive proofs on word length.

If a group is contracting, then the projection map $g\mapsto g_x$ is
not injective, so $x$, as a state of the dual automaton, cannot be
invertible. It is in that sense that contraction and
invertibility-of-dual are opposites.

The following notion is due to Volodymyr
Nekrashevych~\cite{bartholdi-g-n:fractal}. For an automaton group with
states $Q$ and and alphabet $X$, construct the following graph
$\Sigma(\A)$ on the vertex set $X^*$; for all $w\in X^*,x\in X,s\in S$
it has an edge (of the \emph{first} kind) from $w$ to $w^s$, and an
edge (of the \emph{second} kind) from $w$ to $xw$. The edges of the
first kind span the tree $X^*$, and the edges of the second kind span
the disjoint union of the \emph{Schreier graphs} on $X^n$, for all
$n\in\N$.

If $\A$ has an invertible dual, then $\Sigma(\A)$ is a quotient of the
subset of $\Pi$'s Cayley graph spanned by $X^*$. On the other hand, if
$G(\A)$ is contracting, then $\Sigma(\A)$ is Gromov-hyperbolic.

The \emph{limit space} of $G$ is then the hyperbolic boundary of
$\Sigma(\A)$. It can be defined as the equivalence classes of infinite
rays in $\Sigma(\A)$ mutually at bounded distance from each other.

An even stronger property than contraction is that $\sum_{x\in
  X}|g_x|\le\eta|g|+C$ again for some $\eta<1$ and all $g\in G$. Such
a property implies that $G$ has subexponential growth (see
Lemma~\ref{lem:interm}).

A weaker, probabilistic version of this strong contraction property
implies that $G$ is amenable. Namely, if given a uniformly distributed
random group element of length $n$ the distribution of $\sum_{x\in
  X}|g_x|$ has mean less than $\eta n+C$ for some $\eta<1$. There are
groups (for instance $\Gamma$) that satisfy this probabilistic strong
contraction property while have exponential growth; this occurs
because even though there is strong contraction on average, the words
in a geodesic normal form are very far from ``average''.

\subsection{Notation} For $a,b\in G$ and $x,y\in G\cup\Z$ we write
\[a^b=b^{-1}ab;\quad a^{x+y}=a^xa^y;\quad
a^{xy}=(a^x)^y;\quad[a,b]=a^{-1+b}=b^{-a+1}.\]

%%%%%%%%%%%%%%%%%%%%%%%%%%%%%%%%%%%%%%%%%%%%%%%%%%%%%%%%%%%%%%%%
\section{Definitions and Statement of Results}
In this section, $G$ denotes an arbitrary group, and $\Gamma$ denotes
the specific example~\eqref{eq:Gamma} generated by the automaton in
Figure~\ref{fig:autom} (top left).

\subsection{Actions on rooted trees}
Fix a finite alphabet $X=\{\9,\dots,\7\}$. The free monoid $X^*$
naturally has the structure of regular rooted tree, rooted at the
empty word $\emptyset$, with an edge connecting $w$ to $wx$ for all
$w\in X^*$ and $x\in X$. By $wX^*$ we mean the subtree isomorphic to
$X^*$ and attached to the root $\emptyset$ of $X^*$ at its vertex $w$.

Let $W$ denote the automorphism group of $X^*$. Every $g\in W$ induces
a permutation $\pi_g$ of $X$ by restriction, and $g\pi_g^{-1}$ fixes
$X$, so induces for each $x\in X$ an automorphism $g_x$ of $xX^*\cong
X^*$. We therefore have a wreath product decomposition, written
\begin{equation}
  \psi:W\to W\wr\sym X,\quad g\mapsto\pair<g_1,\dots,g_d>\pi_g.
  \label{eq:decomp}
\end{equation}
We will sometimes avoid $\psi$ from the notation for greater clarity.
We also fix a $d$-cycle $(\9,\8,\dots,\7)\in\sym X$.

Let $S$ be a finite subset of $W$; assume that each $s\in S$ appears
exactly once among the $s_x$ for $s\in S,x\in X$; that all other $s_x$
are trivial; and that $\pi_s$ is a power of the fixed $d$-cycle. We
then call the group $G=\langle S\rangle$ a \emph{monomial group}.

As examples on $X=\{\9,\8\}$, we have:
\begin{itemize}
\item $S=\{\lambda^{\pm1},\mu^{\pm1}\}$, with
  \begin{equation}
    \lambda^\psi=\pair<\lambda,1>(\9,\8),\quad\mu^\psi=\pair<\mu^{-1},1>(\9,\8).
    \label{eq:BSV}
  \end{equation}
  This is the ``Brunner-Sidki-Vieira
  group''~\cite{brunner-s-v:nonsolvable}, abbreviated BSV in the
  sequel. It is generated by the automaton in Figure~\ref{fig:autom}
  (top right).
\item $S=\{a,b\}$, with
  \begin{equation}
    a^\psi=\pair<b,1>(\9,\8),\quad b^\psi=\pair<a,1>.
    \label{eq:pink}
  \end{equation}
  This group was discovered by Richard Pink in connection with the
  Galois group of the iterates of the polynomial $z^2-1$. It will be
  called $\Gamma$ in the sequel.
\item More generally, $S=\{a_1,\dots,a_n\}$ with
  \begin{equation}
    a_1^\psi=\pair<a_n,1>(\9,\8),\quad a_i^\psi=\pair<a_{i-1},1>\text{
      or }\pair<a_{i-1},1>.
    \label{eq:mandelbrot}
  \end{equation}
  These groups are the ``iterated monodromy groups'' of polynomials
  $z^2+c$, with $c$ a periodic point in the Mandelbrot set --- see
  Subsection~\ref{subs:coverings}.
\end{itemize}

A group $G$ acting on a rooted tree $X^*$ is \emph{fractal} if for
every $w\in X^*$ the stabilizer of $w$ in $G$ maps to $G$ by
restriction to and identification of $wX^*$ with $wX^*$. The group $G$
is \emph{weakly branch}~\cite{grigorchuk:jibg} if it acts transitively
on $X^n$ for all $n$, and if for each vertex $w\in X^*$ there is a
non-trivial $g\in G$ fixing $X^*\setminus wX^*$ pointwise.

Assume $G$ is finitely generated, and let $|\cdot|$ denote a word
metric on $G$. Then $G$ is \emph{contracting} if there are $\eta<1$
and $C$ such that $|g_x|\le\eta|g|+C$ holds for all $g\in G,x\in X$,
and is \emph{strongly contracting} if there are $\eta<1$ and $C$ such
that $\sum_{x\in X}|g_x|\le\eta|g|+C$ holds for all $g\in G$

\begin{lem}[\cite{bartholdi:upperbd}]\label{lem:interm}
  Let $G$ be strongly contracting, with contraction constant
  $\eta$. Then $G$ has intermediate growth, and its growth function
  $\gamma$ satisfies
  \[\gamma(n)\precsim e^{n^\alpha},\text{ with }\alpha=\frac{\log\#X}{\log(\#X/\eta)}.\]
\end{lem}

By $H^X$ we denote the direct product of $\#X$ copies of $H<W$, acting
independently on the subtrees $xX^*$ for all $x\in X$. Let $G$ be a
fractal group. If it contains a non-trivial subgroup $K$ such that
$K^\psi$ contains $K^X$ in its base group, then $G$ is \emph{weakly
  branch over $K$}; this implies that $G$ is weakly branch.

The main result of this note is the following:
\begin{thm}\label{thm:main}
  Let $G$ be a monomial group. Then $G$ is fractal.  If $G^\psi$ maps
  to a transitive subgroup of $\sym X$, and $\#S\ge2$ with $S$ not of
  the form $\{a,a^{-1}\}$, then $G$ is weakly branch over $G'$.

  If $G'/(G')^X$ is amenable, then $G$ is amenable.
\end{thm}
I do not know whether all monomial groups have exponential growth; the
last two examples do, and this question is open for the BSV
group~\eqref{eq:BSV}.

\subsection{Groups and covering maps}\label{subs:coverings}
The last example above~\eqref{eq:mandelbrot} is a special case of a
construction due to Volodymyr
Nekrashevych~\cite{bartholdi-g-n:fractal}. It was inspired by research
by Richard Pink on Galois groups of iterated polynomials --- see
Point~\eqref{enum:hd} in Theorem~\ref{thm:gamma}.

Let $f$ be a branched self-covering of a Riemann surface $\mathcal S$.
A point $z\in\mathcal S$ is \emph{critical} if $f'(z)=0$, and is a
\emph{ramification point} if it is the $f$-image of a critical point.
The \emph{postcritical set} of $f$ is $\{f^n(z):\,n\ge1,\,z\text{
  critical}\}$.
  
Assume $P$ is finite, and write $\M=\mathcal S\setminus P$. Then $f$
induces by restriction an \'etale map of $\M$.

Let $*$ be a generic point in $\M$, i.e.\ be such that the iterated
inverse $f$-images of $*$ are all distinct. If $f$ has degree $d$,
then for all $n\in\N$ there are $d^n$ points in $f^{-n}(*)$, and all
these points naturally form a $d$-regular tree $U$ with root $*$ and
an edge between $z$ and $f(z)$ for all $z\in U\setminus\{*\}$. Denote
by $X=f^{-1}(*)$ the first level of $U$.

Let $\gamma$ be a loop at $*$ in $\M$. Then for every $v\in U$ at
level $n$ there is a unique lift $\gamma_v$ of $\gamma$ starting at
$v$ such that $f^n(\gamma_v)=\gamma$; and furthermore the endpoint
$v^\gamma$ of $\gamma_v$ also belongs to the $n$th level of $U$.

For any such $\gamma$ the map $v\mapsto v^\gamma$ is a tree
automorphism of $U$, and depends only on the homotopy class of
$\gamma$ in $\pi_1(\M,*)$. We therefore define the \emph{iterated
  monodromy group} $G_U(f)$ of $f$ as the subgroup of $\aut(U)$
generated by all maps $v\to v^\gamma$, as $\gamma$ ranges over
$\pi_1(\M,*)$.

This definition is actually independent of the choice of $*$: if $*'$
is another generic basepoint, with tree $U'$, then choose a path $p$
from $*$ to $*'$. There is then an isomorphism $\phi:U\to U'$ such
that
\[\xymatrix{{\pi_1(\M,*)}\ar[r]^{p_*}\ar[d]_{\text{act}} &
  {\pi_1(\M,*')}\ar[d]^{\text{act}}\\
  {G_U(f)}\ar[r]_{\phi_*} & {G_{U'}(f)}}
\]
commutes; we write $G(f)$ for $G_{U}(f)$, defined up to conjugation in
$\aut(U)$.  Abstractly, $G(f)$ is a presented as a quotient of
$\pi_1(\M,*)$.

We now identify $U$ with the standard tree $X^*$. Enumerate
$X=\{v_1,v_2,\dots,v_d\}$, and choose for each $v\in X$ a path
$\ell_v$ from $*$ to $v$ in $\M$. Consider a loop $\gamma$ at $*$;
then it induces the permutation $v\to v^\gamma$ of $X$, and for each
$v\in X$ its lift $\gamma_v$ at $v$ yields a loop
$\ell_v\gamma_v\ell_{v^\gamma}^{-1}$ at $*$, which depends only on the
class of $\gamma$ in $G(f)$.

We therefore have a natural wreath product decomposition~\eqref{eq:decomp}
\[\phi:G(f)\to G(f)\wr\sym X,\quad g\mapsto \pair<g_1,\dots,g_d>\pi_g,\]
where, if $g$ is represented by a loop $\gamma$ and $v=v_i\in X$, we
have $v^{\pi_g}=v^\gamma$ and $g_i$ is the class of
$\ell_v\gamma_v\ell_{v^\gamma}^{-1}$ in $\pi_1(\M,*)$.

Consider a polynomial self-mapping of the Riemann sphere
$f(z)=z^N+c\in\C[z]$ such that $f^N(0)=0$ for some $N\in\N$. Then
$G(f)$ is a monomial group, acting on $X^*$ where $\#X=d$.  The only
example I consider here is $f(z)=z^2-1$; it gives the group
$G(f)=\Gamma$.

The postcritical set $P$ is $\{0,1,\infty\}$, so $\M$ is a
thrice-punctured sphere and $G$ is $2$-generated.

For convenience, pick as base point $*$ a point close, but not equal,
to $(1-\sqrt5)/2$; then $X=\{x,y\}$ with $x$ close to $*$ and $y$
close to $-*$.

Consider the following representatives of $\pi_1(\M,*)$'s generators:
$a$ is a straight path approaching $-1$, turning a small loop in the
positive orientation around $-1$, and returning to $*$. similarly, $b$
is a straight path approaching $0$, turning around $0$ in the positive
orientation, and returning to $*$.

Let $\ell_x$ be a short arc from $*$ to $x$, and let $\ell_y$ be a
half-circle above the origin from $*$ to $y$.

Let us compute first $f^{-1}(a)$, i.e.\ the path traced by
$\pm\sqrt{z+1}$ as $z$ moves along $a$. Its lift at $x$ moves towards
$0$, passes below it, and continues towards $y$. Its lift at $y$ moves
towards $0$, passes above it, and continues towards $x$. We have
$a_x=b$ and $a_y=1$, so the wreath decomposition of $a$ is
$\phi(a)=\pair<b,1>(1,2)$.

Consider next $f^{-1}(b)$. Its lift at $x$ moves towards $-1$, loops
around $-1$, and returns to $x$. Its lift at $y$ moves towards $1$,
loops and returns to $y$. We have $b_x=a$ and $b_y=1$, so the wreath
decomposition of $b$ is $\phi(b)=\pair<a,1>$.

\begin{figure}
  \begin{center}
    \psfrag{-1}{$-1$}
    \psfrag{0}{$0$}
    \psfrag{a}{$a$}
    \psfrag{b}{$b$}
    \psfrag{x}{$x$}
    \psfrag{y}{$y$}
    \psfrag{a_x}{$a_x$}
    \psfrag{a_y}{$a_y$}
    \psfrag{b_x}{$b_x$}
    \psfrag{b_y}{$b_y$}
    \psfrag{l_x}{$\ell_x$}
    \psfrag{l_y}{$\ell_y$}
    \psfrag{*}{$*$}
    \epsfig{file=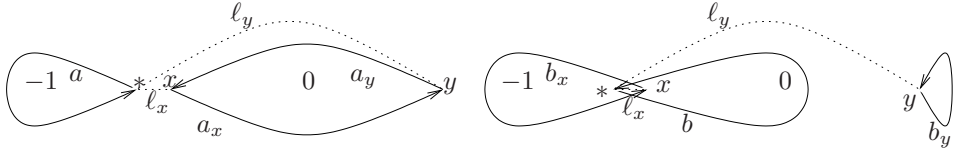,width=5in}
  \end{center}
  \caption{The decomposition of the generators $a,b$ of $\Gamma$.}
  \label{fig:decomposition}
\end{figure}

These paths are presented in Figure~\ref{fig:decomposition}.

\subsection{$\mathbf\Gamma$ and $\mathbf\Delta$}
The next claims concern only the specific example $\Gamma=\langle
a,b\rangle$. For clarity, its action~\eqref{eq:pink} on $\{\9,\8\}^*$ is
given by
\[(\9w)^a =\8w^b,\quad(\8w)^a=\9w,\quad(\9w)^b =\9w^a,\quad(\8w)^b=\8w.\]

\begin{thm}\label{thm:gamma}
  The group $\Gamma$ is
  \begin{enumerate}
  \item fractal, contracting, and weakly branch over
    $\Gamma'$;\label{enum:wb}
  \item torsion free;\label{enum:tf}
  \item of exponential growth, containing $\{a,b\}^*$ as a free
    submonoid;\label{enum:fs}
  \item has as quotients along its lower central series
    $\gamma_1/\gamma_2=\Z^2$, $\gamma_2/\gamma_3=\Z$, and
    $\gamma_3/\gamma_4=\Z/4$. Therefore all successive quotients
    except the first two in the lower central series are finite.
    
    In the lower $2$-central series\footnote{aka ``Zassenhaus
      series'', ``Jennings series'', ``Lazard series'', or ``dimension
      series''} defined by $\Gamma_1=\Gamma$ and
    $\Gamma_{n+1}=[\Gamma,\Gamma_n]\Gamma_{\lfloor n/2\rfloor}^2$, we
    have
    \[\dim_{\F_2}\Gamma_n/\Gamma_{n+1}=\begin{cases}i+2& \text{ if
      }n=2^i\text{ for some }i;\\ 
      \max\setsuch{i+1}{2^i\text{ divides }n}& \text{ otherwise}.
    \end{cases}\]
    \label{enum:cs}
  \item right-orderable, but not bi-orderable\footnote{i.e., there is
      a total order $\le$ on $\Gamma$ with $x\le y\Rightarrow xz\le
      yz$ for all $x,y,z\in\Gamma$, but there is no order satisfying
      $x\le y\Rightarrow wxz\le wyz$ for all
      $w,x,y,z\in\Gamma$.};\label{enum:ro}
  \item not solvable; however, every proper quotient of $\Gamma$ is
    nilpotent-by-(finite 2-group), and every non-trivial normal
    subgroup of $G$ has a subgroup mapping onto $\Gamma$;\label{enum:ns}
  \item has solvable word problem, and is recursively presented as
    \[\Gamma=\langle a,b|\,[[a^p,b^p],b^p],[[b^p,a^{2p}],a^{2p}]\text{
      for all $p$ a power of $2$}\rangle;
    \]
    its Schur multiplier is $H_2(G,\Z)=\Z^\infty$;\label{enum:pr}
  \item Set $f(z)=z^2-1$. Then the closure of $\Gamma$ in the
    profinite group $W$ is the Galois group of $\Lambda/\C(z)$, where
    $\Lambda=\bigcup_{n\ge0}\Lambda_n$ and $\Lambda_n$ is the
    splitting field of $f^n(t)-z$ over $\C(z)$. It has Hausdorff
    dimension\footnote{The automorphism group $W$ of $X^*$ is a
      profinite group, a basis of neighbourhoods of the identity being
      given by the family of pointwise fixators $W_n$ of $X^n$.  For a
      subgroup $G$ of $W$, its \emph{Hausdorff
        dimension}~\cite{barnea-s:hausdorff} is defined by
      \[\dim(G)=\lim_{n\to\infty}\frac{|GW_n/W_n|}{|W/W_n|}.\]}
    $2/3$, and contains the BSV group;\label{enum:hd}
  \item $\Gamma$ has as limit space the Julia set $J$ of $z^2-1$; the
    Schreier graphs of the action of $\Gamma$ on $X^n$ are planar, and
    can be metrized so as to converge to $J$ in the Gromov-Hausdorff
    metric;\label{enum:ls}
  \item The spectrum of the Hecke-type operator\footnote{i.e.\ the
      operator defined as the averaged the sum of the generators of a
      group in a unitary representation. Here $\ell^2(\Gamma)$ denotes
      the left-regular representation of $\Gamma$ by
      left-multiplication on the space of square-summable functions on
      $\Gamma$, and $L^2(X^\omega,\mu)$ denotes the ``natural''
      representation of $\Gamma$ by permutation on the space of
      square-integrable functions on the boundary of the tree $X^*$,
      with $\mu$ the Bernoulli measure.} $\frac14(a+a^{-1}+b+b^{-1})$
    on $L^2(X^\omega,\mu)$ is a Cantor set of null measure, while its
    spectrum on $\ell^2(\Gamma)$ is the interval
    $[-1,1]$.\label{enum:sp}
  \end{enumerate}
\end{thm}

Consider the endomorphism $\sigma:\Gamma\to\Gamma$ given by
$a^\sigma=b,b^\sigma=a^2$, and form the \HNN\ extension
$\Delta=\langle\Gamma,t|\,a^t=b,b^t=a^2\rangle$.
\begin{thm}\label{thm:gammaEG}
  $\Gamma$ and $\Delta$ are amenable, but do not belong to $\BG$.

  $\Gamma$ is infinitely presented, and $\Delta$ has a balanced,
  finite presentation
  \[\Delta = \langle b,t|\,b^{t^2-2},[[[b,t^{-1}],b],b]\rangle.\]
\end{thm}

\subsection{Transitive actions on trees}\label{subs:transitive}
We now consider extension of actions on rooted trees to transitive
actions on regular trees containing the original rooted tree.

\begin{thm}\label{thm:hnn}
  Let $G$ act on $X^*$ and be weakly branch over $K$. Assume that the
  map $K\to K\times1\times\dots\times1$ given by
  $k\mapsto\pair<k,1,\dots,1>$ lifts to an endomorphism
  $\sigma:g\mapsto\pair<g,*,\dots,*>$ of $G$. Then the \HNN\ extension
  $\widetilde G=\langle G,\sigma\rangle$ acts transitively on a
  $(\#X+1)$-regular tree; $\sigma$ is a hyperbolic translation, and
  $G$ is a split quotient of the stabilizer of a vertex $*$ on
  $\sigma$'s axis. Deleting from $*$ the edge on $\sigma$'s axis and
  keeping the connected component of $*$ gives a $\#X$-regular rooted
  tree carrying $G$'s original action.
  
  If furthermore $G$ is contracting, and $G/K$ and $K/K^X$ are both
  finitely presented, then $\widetilde G$ is finitely
  presented\footnote{There is a standard presentation, due to Bass and
    Serre~\cite{serre:trees}, for a group acting transitively on a
    tree. This result should be understood in that spirit.}.
\end{thm}

This result applies to the Pink group $\Gamma$, to the BSV group (see
Equation~\eqref{eq:BSV} or Figure\ref{fig:autom} (top right)), and to
the Grigorchuk group; we start with $G=\Gamma$.

\begin{figure}
  \begin{center}
    \psfrag{xb}{$x^b$}
    \psfrag{yd}{$y^d$}
    \psfrag{w}{$w$}
    \psfrag{x}{$x$}
    \psfrag{y}{$y$}
    \psfrag{z}{$z$}
    \psfrag{bt}{$\tilde b$}
    \psfrag{t}{$t$}
    \psfrag{u}{$u$}
    \epsfig{file=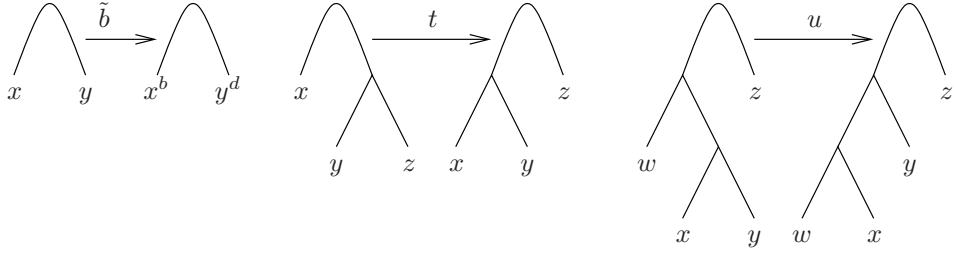}
  \end{center}
  \caption{The action of $\tilde b,t,u$ on the ternary tree, with $e$
    as top edge. $\tilde b$ and $t$ generate the group $\Delta$, while
    $t$ and $u$ generate the Thompson group $F$.}
  \label{fig:thompson}
\end{figure}

Consider a $3$-regular tree $T$; it can be viewed as a rooted binary
tree $U=\{\9,\8\}^*$, in which the root vertex $\emptyset$ was
removed, and its two adjacent edges were replaced by a new edge $e$
joining their extremities $\9,\8$; conversely, a binary tree
isomorphic to $U$ is obtained by inserting a root vertex in the middle
of an edge. Consider the automorphisms $c,d\in W=\aut U$ given by
\[c^\psi=\pair<b,d^2>,\quad d^\psi=\pair<1,c>.\]
The pointwise fixator of $e$ is $W\times W$ acting disjointly on $\9U$
and $\8U$; we still write $\pair<g_1,g_2>$ its elements.  Extend the
action of $\Gamma$ to isometries of $T$ fixing $e$ by letting $a$ act
as $\tilde a=\pair<a,c>$ and letting $b$ act as $\tilde b=\pair<b,d>$.
Note that $[c,d]=1$ so the subgroup of $\aut T$ generated by $\{\tilde
a,\tilde b\}$ is still $\Gamma$.  Let $t$ act by shifting toward the
root in $\8U$ along $\8^\infty$, crossing $e$, and shifting away from
the root in $\9U$ along $\9^\infty$.  In symbols, we have
\[(\8\8w)^t=\8w,\quad(\8\9w)^t=\9\8w,\quad\8^t=\9,\quad(\9w)^t=\9\9w, 
\]
and conjugation by $t$ is given in $\Delta$ by
$\pair<x,\pair<y,z>>^t=\pair<\pair<x,y>,z>$.  It is then easy to
check that $\langle b,t\rangle=\Delta$ in this action, described in
Figure~\ref{fig:pinkaction}. The actions in this setting are given in
Figure~\ref{fig:thompson}.

The Thompson group $F$ is the group of piecewise linear
orientation-preserving self-homeomorphisms of $[0,1]\cap\Z[\frac12]$;
see~\cite{ghys-s:thompson,cannon-f-p:thompson}.  It has a finite,
balanced presentation
\[F=\langle t,u|\,[tu^{-1},u^t],[tu^{-1},u^{t^2}]\rangle.\]
It is known that $F$ is torsion-free, not in the class $\EG$, and does
not contain any non-abelian free subgroup; however, it is open whether
$F$ is amenable.

$[0,1[\cap\Z[\frac12]$ can be identified with $U=\{\9,\8\}^*$ by
mapping the dyadic number $0.x_1\dots x_n$ to $\mathbf{x_1\dots x_n}$.
In this way $F$ acts by homeomorphisms on the boundary of $U$. This
action is described in Figure~\ref{fig:thompson}; note that the
generator $t$ of $F$ acts in the same way as the generator $t$ of
$\Delta$, and $u=\pair<t,1>$ in our notation --- but beware that $u$
is not an isometry of $T$. The arguments in~\cite{rover:fpsg} show
that $\langle\tilde b,t,u\rangle'$ is a finitely presented simple
group. I do not know whether it is amenable, though this question is
probably harder than the corresponding one for $F$.

\begin{figure}
  \begin{center}
    \psfrag{t}{$t$}
    \psfrag{a}{$a$}
    \psfrag{d}{$d$}
    \psfrag{a^d}{$a^d$} 
    \psfrag{0}{\small $0$}
    \psfrag{1}{\small $1$}
    \psfrag{2}{\small $2$}
    \psfrag{3}{\small $3$}
    \psfrag{4}{\small $4$}
    \psfrag{5}{\small $5$}
    \psfrag{6}{\small $6$}
    \psfrag{U}{\large $U$}
    \psfrag{D}{\large $d$}
    \epsfig{file=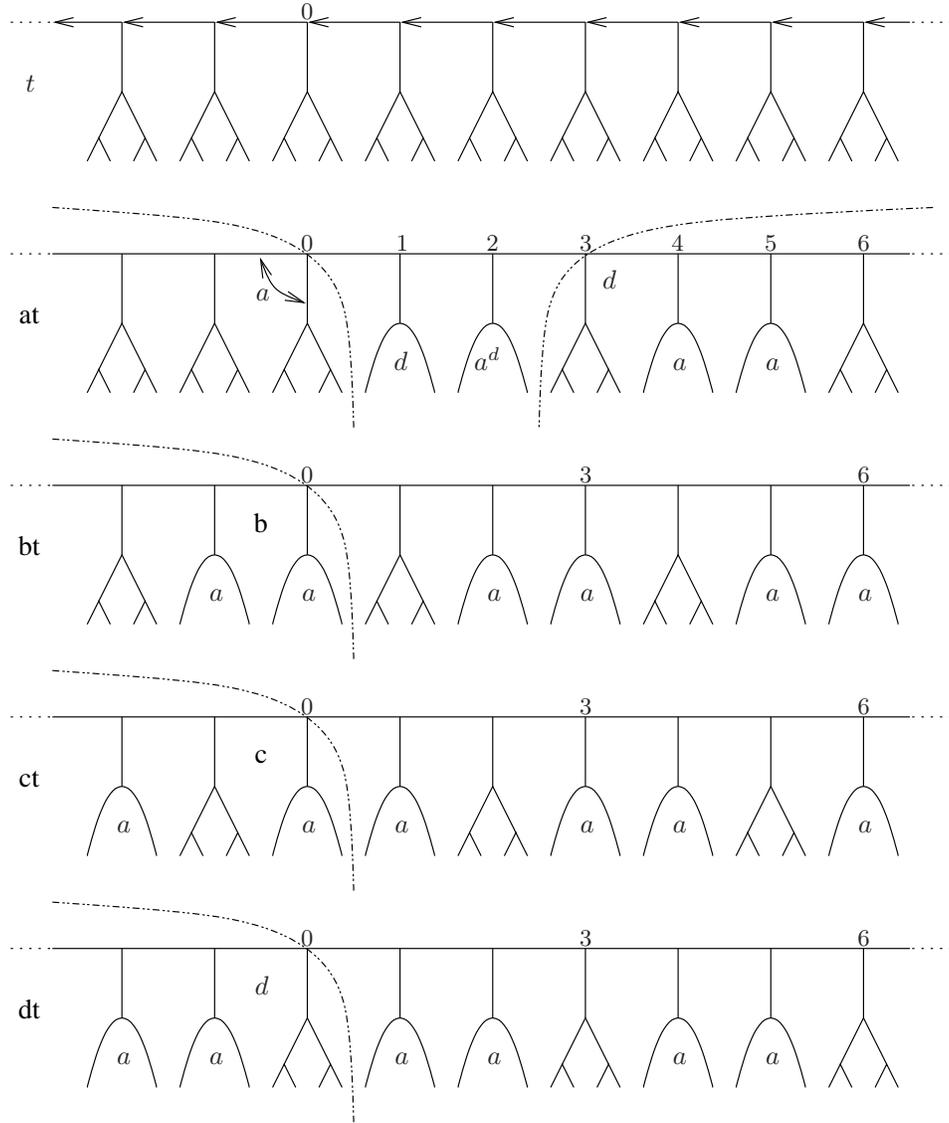,width=5in}
  \end{center}
  \caption{The action of the generators $\tilde a,\tilde b,\tilde
    c,\tilde d,t$ of $\widetilde G$ on the $3$-regular tree $T$.}
  \label{fig:grigaction}
\end{figure}

Consider next the Grigorchuk group $G$ from Figure~\ref{fig:autom}
(bottom left). It may be defined as $G=\langle a,b,c,d\rangle$ acting
on $\{\9,\8\}^*$, with
\begin{equation}
  a^\psi=\pair<1,1>(\9,\8),\quad b^\psi=\pair<a,c>,\quad
  c^\psi=\pair<a,d>,\quad d^\psi=\pair<1,b>.
  \label{eq:G}
\end{equation}
This group is contracting, and even strongly
contracting~\cite{bartholdi:upperbd}, with $|g_1|+|g_2|\le\eta(|g|+1)$
for $\eta\cong0.811$ the real root of $X^3+X^2+X-2$. It is therefore
of intermediate growth, of rate at most $e^{n^{0.768}}$, and hence is
amenable. It embeds in the finitely presented group $\widetilde
G=\langle G,t\rangle$, with $t$ acting by conjugation as the
endomorphism $\sigma:G\to G$ given by
\[a^\sigma=c^a,\quad b^\sigma=d,\quad c^\sigma=b,\quad d^\sigma=c.\]
Consider the following isometries of the $3$-regular tree $T$
described above:
\[\tilde a=\pair<a,\pair<d,\pair<a^d,d>>>,\quad\tilde b=\pair<b,d>,\quad
\tilde c=\pair<c,c>,\quad\tilde d=\pair<d,b>.
\]
Then $G\cong\langle\tilde a,\tilde b,\tilde c,\tilde d\rangle$, and
$\widetilde G$ is generated by $t\tilde a,\tilde b,\tilde c,\tilde d$
and a hyperbolic element $t$ moving toward the root along $\8^\infty$
in $\8U$ and away from the root along $\8^\infty$ in $\9U$; we have
$\pair<x,\pair<y,z>>^t=\pair<\pair<y,x>,z>$. The action of $\widetilde
G$ is described in Figure~\ref{fig:grigaction}.

A presentation of $\widetilde G$ with $2$ generators and $4$ relators,
obtained using $c=a^{ta},b=a^{tat},d=a^{tat^2}$, is
\[\widetilde G=\langle a,t|\,a^2,a^{tat^2+tat+ta},a^{(1+ta)8},
a^{(1+tat^2+(1+ta)2)4}\rangle.\]

\subsection{Reddite Caesare}
Some of the results in Theorem~\ref{thm:gamma} were obtained
independently by Rostislav Grigorchuk and Andrzej \.Zuk, whom the
author thanks for their communication. The proof technique follows
ideas appearing in the original works of Rostislav Grigorchuk, Said
Sidki~\cite{brunner-s-v:nonsolvable} and Edmeia da
Silva~\cite{silva:phd}. The author is also extremely grateful to
Professors de la Harpe, Grigorchuk and Nekrashevych for their generous
sharing of knowledge and ideas.

%%%%%%%%%%%%%%%%%%%%%%%%%%%%%%%%%%%%%%%%%%%%%%%%%%%%%%%%%%%%%%%%
\section{Proofs}
We use $S$ as a natural generating set of $G$, and write $|w|$ the
length of a word, and $|g|$ the minimal length of a group element.
Most of the proofs follow by induction on $|g|$.

\begin{proof}[Proof of Theorem~\ref{thm:main}]
  Write $K=G'$. For $s,t\in S$, pick $s',t'\in S$ such that
  $s'_x=s,t'_y=t$, and let $n$ be the order of $\pi_t$. Then
  $[s',(t')^n]^\psi$ will have precisely one non-trivial coordinate,
  containing $[s,t]$. By conjugating, $K$ contains $K^X$. Finally
  $K\neq1$ by our assumption.
  
  Consider next the set $F_n$ of freely reduced words of length $n$
  over $S$, and the subset $N_n$ of words evaluating to $1$ in $G$.
  $F=\bigcup_{n\ge0}F_n$ is the free group on $S$, and
  $N=\bigcup_{n\ge0}N_n$ is the kernel of the natural map $F\to G$.
  
  \begin{lem}[Kesten~\cite{kesten:rwalks}; Grigorchuk~\cite{grigorchuk:rw}]
    $G$ is amenable if and only if $\#N_n/\#F_n>\rho^n$ for all
    $\rho<1$ and all $n$ even and large enough.
  \end{lem}

  The decomposition map $\psi:G\to G\wr\langle(\9,\dots,\7)\rangle$
  induces a map $F\to F\wr\langle(\9,\dots,\7)\rangle$ on freely reduced
  words, again written $w^\psi=\pair<w_1,\dots,w_d>\pi_w$.  By
  construction, we have $|w_1|+\dots+|w_d|\le|w|$; and usually the
  inequality is strict: since $G$ is weakly branch, there are
  non-trivial reduced words $u,v$ with $u_x=1$ for all $x\neq1$ and
  $v_1=1$; then $w=[u,v]$ has positive length but $w_1=\dots=w_d=1$.
  
  The cancellation that occurs in the $w_d$ is determined by the
  \emph{local} rules specifying the decomposition of generators.
  Therefore, if $w$ is chosen uniformly at random in $F_n$ with $n$
  large, then $w_x$ will again be uniformly distibuted within
  $F_{|w_x|}$, and the length of each $w_x$ will follow a binomial
  distribution; hence $|w_1|+\dots+|w_d|$ will also follow a binomial
  distribution.
  
  Assume that the mean of this distribution is $\mu n$ and its
  variance is $\mu(1-\mu)n/\eta$. This means that the probability that
  $w$ of length $n$ yields via $\psi$ freely reduced words
  $w_1,\dots,w_d$ of total length $m$ is
  \begin{equation}
    C_{m,n}^{\eta,\mu}=\eta\binom{\eta n}{\eta m}\mu^{\eta m}(1-\mu)^{\eta(n-m)}.
    \label{eq:C}
  \end{equation}
  By the above argument we have $\mu<1$, although the precise value is
  unimportant for the present.
  
  For a subgroup $A<G$ we write $p_A(n)$ the probability that $w$
  chosen uniformly at random in $F_n$ evaluates to $1\in G$,
  conditionally on knowing that that it evaluates to an element of
  $A$. For $A<B<G$ we write $p_{A/B}(n)=p_A(n)/p_B(n)$.
  
  Assume now for contradiction that $p(n)$ decays exponentially at
  rate $\rho$, say $(\rho-\epsilon)^n<p(n)<(\rho+\epsilon)^n$ for any
  $\epsilon>0$, provided $n$ is large enough. Since $G/K$ is abelian,
  $p_{G/K}(n)$ decays subexponentially, so we also have
  $(\rho-\epsilon)^n<p_K(n)<(\rho+\epsilon)^n$. Then for large $n$
  \begin{align*}
    (\rho+\epsilon)^n&\ge p(n)\ge p_{G/K}(n)p_K(n),\\
    \text{and }p_K(n)&\ge p_{K/K^X}(n)\sum_{\substack{0\le m\le n\\
        i_1+\dots+i_d=m}}C_{n,m}^{\eta,\mu}p_K(i_1)\dots p_K(i_d);\\
    \intertext{Writing $E(n)$ a function that decays subexponentially,
      and takes into account both the $\approx\binom md$ ways of
      partitioning $m$ in $d$ parts, and $p_{K/K^X}(n)$, of
      subexponential decay since $K/K^X$ is assumed to be amenable,}
    p_K(n)&\ge E(n)\eta\sum_{m=0}^n\binom{\eta n}{\eta
      m}\mu^{\eta m}(1-\mu)^{\eta(n-m)}(\rho-\epsilon)^m\\
    &\approx E(n)\sum_{m'=0}^{\eta n}\binom{\eta
      n}{m'}(\mu\sqrt[\eta]{\rho-\epsilon})^{m'}(1-\mu)^{\eta n-m'}\\
    &=E(n)\left((1-\mu)+\mu\sqrt[\eta]{\rho-\epsilon}\right)^{\eta n}.
  \end{align*}
  Letting next $n$ tend to $\infty$ and taking $n$th roots, we get
  $\sqrt[\eta]{\rho+\epsilon}\ge(1-\mu)+\mu\sqrt[\eta]{\rho-\epsilon}$
  and hence $\rho\ge1$, since $\epsilon>0$ is arbitrary, $\mu<1$, and
  $\eta>0$.
\end{proof}
We note that the parameters $\mu,\eta$ were experimentally found to be
$\mu\approx0.699$ and $\eta\approx0.326$ for the Pink group
$\Gamma$, and $\mu\approx0.781$ and $\eta\approx0.282$ for the BSV
group. These values were obtained by a Monte-Carlo simulation using
$1\,000\,000$ words of length $50\,000$.

I now proceed with the proof of Theorem~\ref{thm:gamma} describing
algebraic properties of $\Gamma$. Alternate proofs of some of the
points were found independently by Grigorchuk and \.Zuk, and appear
in~\cite{grigorchuk-z:torsionfree}.

\begin{figure}
  \begin{center}
    \psfrag{G1}[c][c]{$\gf_1$}
    \psfrag{G2}[c][c]{$\gf_2$}
    \psfrag{G3}[c][c]{$\gf_3$}
    \psfrag{G4}[c][c]{$\gf_4$}
    \psfrag{G5}[c][c]{$\gf_5$}
    \psfrag{G6}[c][c]{$\gf_6$}
    \psfrag{A1}[c][c]{$A_1$}
    \psfrag{A2}[c][c]{$A_2$}
    \psfrag{A3}[c][c]{$A_3$}
    \psfrag{A4}[c][c]{$A_4$}
    \psfrag{A5}[c][c]{$A_5$}
    \psfrag{A6}[c][c]{$A_6$}
    \psfrag{B1}[c][c]{$B_1$}
    \psfrag{B2}[c][c]{$B_2$}
    \psfrag{B3}[c][c]{$B_3$}
    \psfrag{B4}[c][c]{$B_4$}
    \psfrag{B5}[c][c]{$B_5$}
    \psfrag{B6}[c][c]{$B_6$}
    \epsfig{file=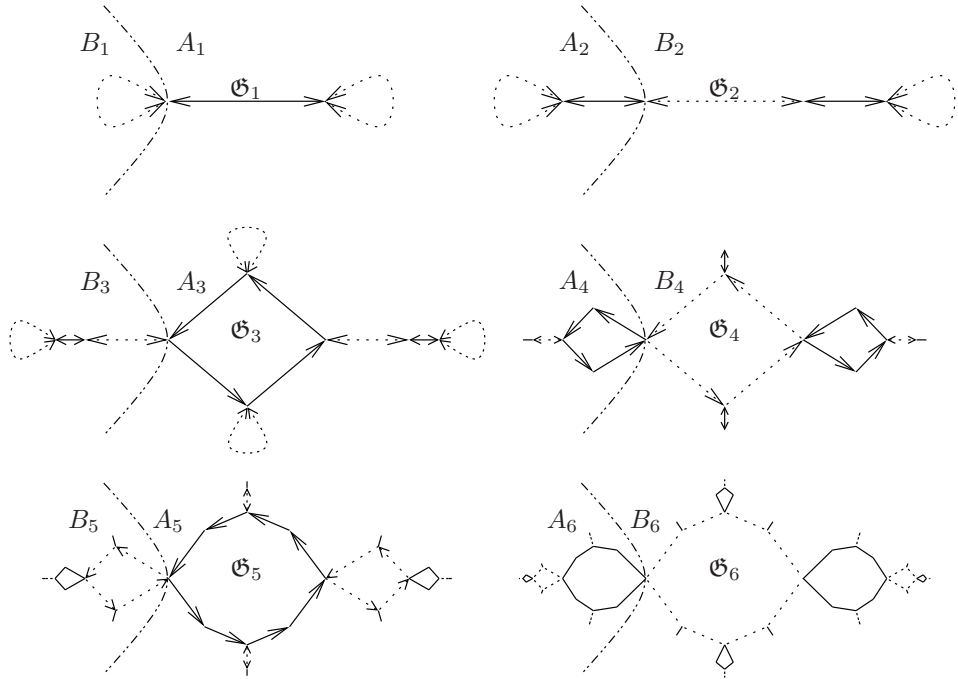}
  \end{center}
  \caption{The Schreier graphs $\gf_n$ for $1\le n\le 6$. The solid
    lines represent $a$'s action on $X^n$, and the dotted lines
    represent $b$'s action. All vertices have degree $4$; the $b$
    loops are represented only for $n\le3$.}
  \label{fig:schreier}
\end{figure}

\begin{figure}
  \begin{center}
    \epsfig{file=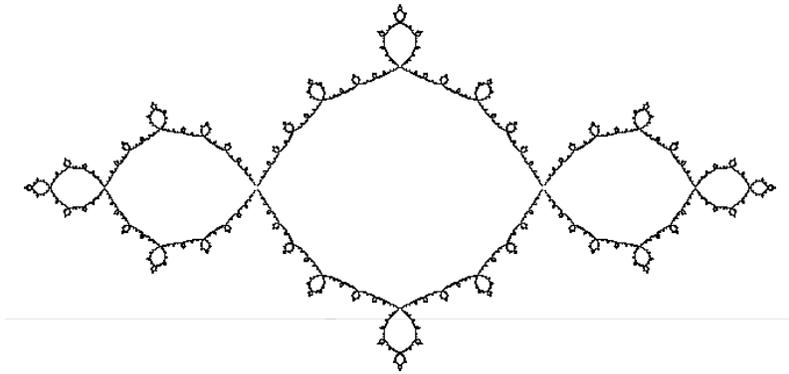,width=100mm}
  \end{center}
  \caption{The Julia set of the polynomial $z^2-1$.}
  \label{fig:julia}
\end{figure}

\begin{figure}
  \psfrag{\2610.6}[tc]{$-0.6$}
  \psfrag{\2610.4}[tc][tc]{$-0.4$}
  \psfrag{\2610.2}[tc][tc]{$-0.2$}
  \psfrag{0}[tc][bc]{$0$}
  \psfrag{0.2}[tc][Bc]{$0.2$}
  \psfrag{0.4}[tc][cc]{$0.4$}
  \psfrag{0.6}[tc][bc]{$0.6$}
  \psfrag{0.8}[tc][bc]{$0.8$}
  \psfrag{1}[tc]{$1$}
  \begin{center}
    \epsfig{file=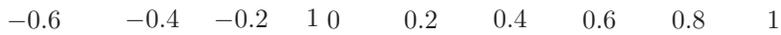,width=300pt,height=20pt}
  \end{center}
  \caption{The spectrum of $\pi$, in its level-$6$ approximation.}
  \label{fig:spectrum}
\end{figure}

For convenience, we write $c=[a,b]$, $d=[c,a]$ and $e=[d,a]$ in
$\Gamma$.
\begin{proof}[Point~\eqref{enum:wb} of Theorem~\ref{thm:gamma}]
  $\Gamma$ is fractal and weakly branch by Theorem~\ref{thm:main}.  By
  letting $a$ have length $1$ and $b$ have length $\sqrt2$, we have
  $|g_x|\le(|g|+1)/\sqrt2$ for all $g\in G,x\in X$; hence $\Gamma$ is
  contracting.
  
  Next, we prove by induction on length of words that we have
  $\Gamma/\Gamma'\cong\Z^2$ generated by the images of $a,b$; and
  $\Gamma'/(\Gamma'\times\Gamma')\cong\Z$, generated by the image of
  $c$.
  
  Assume for contradiction that $a^mb^n\in\Gamma'$ with $|m|+|n|$
  minimal.  Then clearly $m$ is even, say $m=2p$. We have, for some
  $k\in\Z$ with $|k|\le|m|+|n|$,
  \[a^mb^n=\pair<b^pa^n,b^p>=\pair<g,h>c^k=\pair<ga^k,ha^{-bk}>\]
  and therefore $a^{n-k}b^p$ and $a^kb^p$ both belong to $\Gamma'$. This
  contradicts our assumption on minimality.
  
  Assume next that $c^k\in \Gamma'\times\Gamma'$ with $|k|$ minimal.
  Then $a^k\in\Gamma'$ which contradicts the second claim.
\end{proof}
\begin{proof}[Point~\eqref{enum:tf}]
  Since $\Gamma$ acts on the binary tree, it is residually a $2$-group,
  and its only torsion must be $2$-torsion. Assume for contraction
  that $\Gamma$ contains an element $g$ of order $2$, of minimal norm.
  
  By the previous point, $a$ and $b$ are of infinite order.  We may
  therefore assume $|g|\ge2$.  If $g$ fixes $\9$, then its
  restrictions $g_x,x\in X$ are shorter, and at least one of them has
  order $2$, contradicting $|g|$'s minimality.
  
  If $g$ does not fix $x$, then we may write $g=\pair<g_1,g_2>a$ for
  some $g_1,g_2\in\Gamma$. We then have $h=g^2=\pair<g_1
  bg_2,g_2g_1b>=1$ and therefore $g_2g_1b=1$. Now for any element $h$
  fixing $\9$ we have $h_1h_2\in\langle a,b^2,\Gamma'\rangle$; this
  last subgroup does not contain $b$ by the previous Point, so we have
  a contradiction.
\end{proof}
\begin{proof}[Point~\eqref{enum:fs}]
  Consider two words $u,v$ in $\{a,b\}^*$ that are equal in $\Gamma$,
  and assume $|u|+|v|$ is minimal. We have $u_1=v_1$ and $u_2=v_2$ in
  $\Gamma$, which are shorter relations, so we may assume these words
  are equal by induction.
  
  Now if $u_1$ and $v_1$ start with the same letter $a$ or $b$, this
  implies that $u$ and $v$ also start with the same letter $b$ or $a$
  respectively, and cancelling these letters would give a shorter pair
  of words $u,v$ equal in $\Gamma$.
  
  It follows that $\{a,b\}^*$ is a free submonoid, and hence that
  $\Gamma$ has exponential growth.
\end{proof}  
\begin{proof}[Point~\eqref{enum:cs}]
  This follows from writing generators for $\gamma_i$, and using
  induction on length. Writing $c=[a,b]$, $d=[c,a]$ and $e=[d,a]$,
  \begin{align*}
    \gamma_1 &= \langle a,\,b\rangle;\\
    \gamma_2 &= \langle c=[a,b]=(a,a^{-b}),\,
    c^{-1-a}=\pair<c,1>,\,c^{-a^{-1}-1}=\pair<1,c>\rangle;\\
    \gamma_3 &= \langle d=[c,a],\,e=[d,a],\,[e^{-1},b]=\pair<d,1>,\,\pair<e,1>,\,\pair<1,d>,\,\pair<1,e>\rangle;\\
    \gamma_4 &= \langle d^4,\,e,\pair<d,1>,\,\pair<e,1>,\,\pair<1,d>,\,\pair<1,e>\rangle.
  \end{align*}
  Only $d^4\in\gamma_4$ deserves some justification; writing $\equiv$
  for congruence modulo $\gamma_4$, we have
  \begin{align*}
    d^2 &\equiv d^2e=b^{-1}a^{-1}ba^{-2}b^{-1}aba^2\\
    &=b^{-1}aba^{-2}b^{-1}a^{-1}ba^2\text{ using the relation }[a^{2b},a^2]=1\\
    &=(d^2e)^{-a^{-b}}\equiv d^{-2}.
  \end{align*}
  For the $2$-central series see~\cite{bartholdi:bsvlcs}, where the
  same answer is proven for the BSV group.
\end{proof}
\begin{proof}[Point~\eqref{enum:ro}]
  We consider for all $n\in\N$ the subgroups
  $\Gamma_n=(\Gamma')^{X^n}$ of $\Gamma$. Then
  $\Gamma/\Gamma_0\cong\Z^2$ and $\Gamma_n/\Gamma_{n+1}\cong\Z^{2^n}$
  are both right-orderable, and $\bigcap_{n\ge0}\Gamma_n=1$. Define a
  right order on $\Gamma$ by
  \[x\le y\Leftrightarrow x=y\text{ or }xy^{-1}<1\text{ in
  }\Gamma_n/\Gamma_{n+1},\text{ where }n\text{ is maximal with
  }xy^{-1}\in\Gamma_n.
  \]
  Note that this is not a bi-ordering, since $\Gamma_{n+1}$ is not
  central in $\Gamma_n$. That no bi-ordering exists follows from
  $d^2e$ being conjugate to $(d^2e)^{-1}$, see Point~\eqref{enum:cs}.
\end{proof}
\begin{proof}[Point~\eqref{enum:ns}] 
  The calculations in Point~\eqref{enum:cs}] show that $\Gamma''$ is
  the normal closure of $[c,\pair<c,1>]=\pair<d,1>$; therefore
  $\Gamma''=\gamma_3\times\gamma_3$, and so $G'''>G''\times G''$;
  hence $G^{(n)}>G^{(n-1)}\times G^{(n-1)}$ for all $n$. Assume for
  contradiction that $G$ is solvable; this means $G^{(n)}=1$ for some
  minimal $n$, a contradiction with the above statement.
  
  Now consider a non-trivial normal subgroup $N$ of $G$.
  By~\cite[Theorem~4]{grigorchuk:jibg}, we have $(\gamma_3)^{X^n}<N$
  for some $n$. Since $G$ is an abelian-by-(finite $2$) extension of
  $(\gamma_2)^{X^n}$, we conclude that $G/N$ is nilpotent-by-(finite
  $2$).
  
  To show that a subgroup of $N$ maps onto $\Gamma$, since $N$
  contains $(\gamma_3)^{X^n}$, it is sufficient to show that
  $\gamma_3$ maps onto $\Gamma$. Now $\gamma_3$ is the normal closure
  of $d=[[a,b],a]$ in $\Gamma$. We have
  $d^\psi=\pair<a^{-1-b},a^{2b}>$, and
  $a^{(-1-b)\psi}=\pair<b^{-1}a,b^{-1}a^{-1}>$ and
  $a^{2b\psi}=\pair<b^a,b>$; therefore projection twice on the first
  factor maps $\gamma_3$ to $\Gamma$.
\end{proof}
\begin{proof}[Point~\eqref{enum:pr}]
  Let $F$ be the free group on $\{a,b\}$, and write $\Gamma=F/R$.
  Then~\eqref{eq:pink} defines a homomorphism $F\to F\wr\sym X$.
  Letting $\sigma$ denote the $F$-endomorphism $a\mapsto b,b\mapsto
  a^2$ we have a diagram
  \[\xymatrix{{\langle a^2,b,b^a\rangle}\ar[r]^{\pi}\ar@{-}[d] & {F\times F}\ar@{-}[d]\\
    {R}\ar[r]\ar@{-}[d] & {R\times R}\ar@{-}[d]\\
    {R^\sigma}\ar[r]^{\cong} & {R\times1}}
  \]
  Set $R_0=1<F$ and inductively $R_{n+1}=(R_n\times R_n)^{\pi^{-1}}$.
  Then $R=\bigcup_{n\ge0}R_n$ because $\Gamma$ is contracting, and we
  have $R_{n+1}=(R_1R_n^\sigma)^F$. Since
  $R_1=\langle[b,b^{a^i}]:\text{ odd }i\rangle^F$, we have
  $R=\langle[b^{a^i},b]^{\sigma^n}:n\in\N,\text{ odd }i\rangle$. Now
  $b^{a^i}\equiv[b^{-1},a^{-2}]^ab^{a^{i-2}}$ using the relation
  $[a^{2b},a^2]=[b^a,b]^\sigma$; therefore $[b^{a^i},b]$ follows from
  $[b^{a^{i-2}},b]$ and $[[b^{-1},a^{-2}]^a,b]$, which in turn is a
  consequence of $[b^a,b]$; the presentation of $\Gamma$ follows.

  The Schur multiplier of $G$ is $(R\cap[F,F])/[R,F]$, by Hopf's
  formula. Writing $R=\langle[b,b^a]\rangle R^\sigma$, we get
  \[R/[R,F]=\langle[b,b^a]\rangle\oplus R^\sigma/[R,F]^\sigma,\]
  so $H^2(G,\Z)\cong\Z^\infty$ with $\sigma$ acting on it as a
  one-sided shift.
\end{proof}
\begin{proof}[Point~\eqref{enum:hd}]
  Write $Q_n$ the quotient $\Gamma W_n/W_n$. Induction shows that $a$
  has order $2^{\lceil n/2\rceil}$ in $Q_n$, and $b$ and $[a,b]$ have
  order $2^{\lfloor n/2\rfloor}$ in $Q_n$; hence $Q_n'$ has index
  $2^n$ in $Q_n$ and $Q_{n-1}'\times Q_{n-1}'$ has index $2^{\lfloor
    n/2\rfloor}$ in $Q_n'$. Since $|Q_0|=1$, we get
  \[|Q_n|=2^{\frac23(2^n+\frac12\lfloor\frac{3n}2\rfloor-1)};\]
  since $|W/W_n|=2^{2^n-1}$, we have $\dim(G)=\frac23$.
  
  The generator $\mu$ of the BSV group is $b^{-1}a\in\Gamma$, since
  $b^{-1}a=\pair<a^{-1}b,1>(\9,\8)$ satisfies $\mu$'s
  recursion~\eqref{eq:BSV}. We do not have $\tau\in\Gamma$; but
  defining $c_n\in\Gamma$ by the recursion $c_0=1$ and
  $c_{n+1}=\pair<1,[b,a]c_n>$, we have
  $bac_{n+1}=\pair<ab[b,a]c_n,1>(\9,\8)$; hence setting
  $c=\lim_{n\to\infty}c_n\in\overline\Gamma$ we have
  $bac=\pair<bac,1>(\9,\8)=\tau\in\overline\Gamma$.
\end{proof}
\begin{proof}[Point~\eqref{enum:ls}]
  Consider next the Schreier graphs $\gf_n$ with $X^n$ as the vertex
  set.  $\gf_n$ is constructed as follows: it is built of two parts
  $A_n,B_n$ connected at a distinguished vertex. Each of these parts
  is $4$-regular, except at the connection vertex where each is
  $2$-regular, and $A_n$ contains only the $a^{\pm1}$-edges while
  $B_n$ contains only the $b^{\pm1}$-edges.
  
  $A_0$ and $B_0$ are the graphs on $1$ vertex with a single loop of
  the appropriate label.
  
  If $n=2k$ is even, then $B_{2k+1}=B_{2k}$, and $A_{2k+1}$ is
  obtained by taking an $a$-labelled $2^{k+1}$-gon
  $v_0,\dots,v_{2^{k+1}-1}$, and attaching to each $v_i$ with $i\neq0$
  a copy of $B_{2j}$ where $2^j||i$. Its distinguished vertex is
  $v_0$.
  
  If $n=2k-1$ is odd, then $A_{2k}=A_{2k-1}$, and $B_{2k}$ is obtained
  by taking a $b$-labelled $2^k$-gon $v_0,\dots,v_{2^k-1}$, and
  attaching to each $v_i$ with $i\neq0$ a copy of $A_{2j+1}$ where
  $2^j||i$. Its distinguished vertex is $v_0$.
  
  The first Schreier graphs $\gf_n$ of $G$ are drawn in
  Figure~\ref{fig:schreier}. Compare with the Julia set in
  Figure~\ref{fig:julia}.
\end{proof}
\begin{proof}[Point~\eqref{enum:sp}]
  Consider first the spectrum on $\hilb=L^2(X^\omega,\mu)$. Since
  $X^\omega=\9X^\omega\sqcup\8X^\omega$, we may decompose $a,b$ and
  write them as $2\times2$-matrices over $B(\hilb)$. We have
  $a=(\begin{smallmatrix}0&b\\1&0\end{smallmatrix})$ and
  $b=(\begin{smallmatrix}a&0\\0&1\end{smallmatrix})$; finite
  approximations $a_n,b_n$ can be obtained by expanding to
  $2^n\times2^n$-matrices and replacing all $a$'s and $b$'s by $1$; we
  have
  \[a_0=b_0=(1),\quad a_{n+1}=\begin{pmatrix}0&b_n\\1&0\end{pmatrix},
  \quad b_{n+1}=\begin{pmatrix}a_n&0\\0&1\end{pmatrix}.
  \]
  Introduce for $n\ge0$ the following homogeneous polynomials of
  degree $2^n$:
  \[Q_n(\lambda,\mu,\nu)=\det\left(\lambda+\mu(a_n+a_n^{-1})+\nu(b_n+b_n^{-1})\right).\]
  Then the solution of $Q_n(\lambda,-\frac14,-\frac14)=0$ is the
  spectrum of the Hecke-type operator
  $\frac14(a_n+a_n^{-1}+b_n+b_n^{-1})$ of $\Gamma$'s action on $\C
  X^n$; this is also the spectrum of the Schreier graph $\gf_n$.

  Define the homogeneous polynomial mapping $F:\R^3\to\R^3$ by
  \[(\lambda,\mu,\nu)\mapsto(\lambda^2+2\lambda\nu-2\mu^2,\lambda\nu+2\nu^2),-\mu^2).\]
  Then $Q_n$ is given, for $n\ge1$, by
  \begin{align*}
    Q_0(\lambda,\mu,\nu)&=\lambda+2\mu+2\nu;\\
    Q_1(\lambda,\mu,\nu)&=Q_0(\lambda,\mu,\nu)\cdot(\lambda-2\mu+2\nu);\\
    Q_{n+1}(\lambda,\mu,\nu)&=\det\begin{pmatrix}\lambda+\nu(a_n+a_n^{-1})
      & \mu(1+b_n)\\ \mu(b_n^{-1}+1) & \lambda+2\nu\end{pmatrix}\\
    &= \det\left((\lambda+\nu(a_n+a_n^{-1}))(\lambda+2\nu)
      -\mu^2(1+b_n)(b_n^{-1}+1)\right)\\
    &= Q_n(\lambda^2+2\lambda\nu-2\mu^2,\lambda\nu+2\nu^2),-\mu^2)\\
    &= Q_n(F(\lambda,\mu,\nu)).
  \end{align*}
  
  Define $K$ as the closure of the set of all backwards $F$-iterates
  of $\{Q_1=0\}$. Then the spectrum of $\pi$ is the intersection of
  the line $\{\mu=\nu=-\frac14\}$ with $K$, and is easily seen to be a
  Cantor set --- see Figure~\ref{fig:spectrum}.

  By contrast, the spectrum of $\ell^2(\Gamma)$ is the interval
  $[-1,1]$ by~\cite{higson-k:bc}, since $\Gamma$ is amenable and
  torsion-free (and hence satisfies the Baum-Connes conjecture).
\end{proof}

My proof that $\Gamma$ is not in $\BG$ is inspired
by~\cite{grigorchuk-z:torsionfree}. Define for ordinals $\alpha$ the
following subclasses of $\BG$: first, $\BG_0$ is the class of groups
locally of subexponential growth. Let $\BG_{\alpha+1}$ be the class of
subgroups, quotients, extensions and direct limits of groups in
$\BG_\alpha$, and for a limit ordinal $\beta$ set
$\BG_\beta=\bigcup_{\alpha<\beta}\BG_\alpha$. Note that it is actually
not necessary to consider subgroups and quotients in the inductive
construction of $\BG_{\alpha+1}$.

\begin{proof}[Proof of Theorem~\ref{thm:gammaEG}]
  $\Gamma$ is amenable by Theorem~\ref{thm:main}.  Since $\Delta$ is
  an ascending extension of $\Gamma$, it is also amenable.
  
  Assume $\Gamma\in\BG$ for contradiction. Then $\Gamma\in\BG_\alpha$
  for some minimal ordinal $\alpha$, which of course is not a limit
  ordinal. Since $\Gamma$ has exponential growth, we have $\alpha>0$.
  
  By minimality of $\alpha$, $\Gamma$ cannot be a subgroup or quotient
  of a group in $\BG_{\alpha-1}$. It cannot be a direct limit, since
  it is finitely generated. Therefore there are $N,Q\in\BG_{\alpha-1}$
  with $G/N=Q$. Now by Theorem~\ref{thm:gamma}, Point~\eqref{enum:ns},
  $N$ has admits a subgroup mapping onto $\Gamma$, so
  $\Gamma\in\BG_{\alpha-1}$, a contradiction.
  
  The presentation of $\Delta$ is obtained as the \HNN\ extension of
  $\Gamma$ identifying $\Gamma$ and $\Gamma^\sigma$. To $\Gamma$'s
  presentation we add a generator $t$ and relations $a^t=b,b^t=a^2$;
  and note then that of the relations of $\Gamma$ all can be removed
  but the first, and $a$ can be removed from the generating set and
  replaced by $b^{t^{-1}}$ in $[b^a,b]$ and $b^ta^{-2}$.
\end{proof}

\begin{proof}[Proof of Theorem~\ref{thm:hnn}]
  Consider the sequence of trees with basepoint $T_n=(X^*,\9^n)$, where
  $\9^n$ is the leftmost vertex at level $n$ in $X^*$. The direct limit
  $\injlim T_n$ is a $(\#X+1)$-regular tree with a distinguished vertex
  $*$. The tree injection $T_n\to T_{n+1}$ given by $w\mapsto 1w$
  extends to an invertible hyperbolic isometry $t$ of $T$.
  
  Let $g\in G$ act on $T_n$ as $g^\sigma$ acts on $X^*$; this action
  extends to the limit $T$, and we have $g^\sigma=g^t$, so $\langle
  G,t\rangle$ is an \HNN\ extension.
  
  Let $U$ denote the connected component of $T\setminus\{\text{the edge
  }\9\text{ at }*\}$. Then $U$ is naturally isomorphic to $X^*$ and
  carries the original action of $G$. Therefore restriction to $U$
  gives a split epimorphism from the stabilizer of $U$ to $G$.

  Any $v\in T$ can be mapped to a vertex in $U$ by a sufficiently
  large power of $t$; then since $G$ is transitive on $X^n$ for all
  $n$, it further can be mapped to some vertex $\9^n$; and mapped to
  $*$ by $t^{-n}$; therefore $\widetilde G$ acts transitively on $T$.
  
  Finally, if $G/K$ and $K/K^X$ are finitely presented and $G$ is
  contracting, then $\widetilde G$ is finitely presented, by the
  argument in~\cite{bartholdi:lpres}.
\end{proof}  

\bibliography{mrabbrev,people,math,grigorchuk,bartholdi}

\def\nop#1{}\font\cyr=wncyr8\def\cprime{$'$}
\providecommand{\bysame}{\leavevmode\hbox to3em{\hrulefill}\thinspace}
\begin{thebibliography}{CGH99}

\bibitem[Ady82]{adyan:rw}
Sergei~I. Adyan, \emph{Random walks on free periodic groups}, Izv. Akad. Nauk
  SSSR Ser. Mat. \textbf{46} (1982), no.~6, 1139--1149, 1343.

\bibitem[Ale83]{aleshin:free}
Sergei~V. Ale{\v{s}}in, \emph{A free group of finite automata}, Vestnik Moskov.
  Univ. Ser. I Mat. Mekh. (1983), no.~4, 12--14.

\bibitem[Bar98]{bartholdi:upperbd}
Laurent Bartholdi, \emph{The growth of {Grigorchuk}'s torsion group}, Internat.
  Math. Res. Notices \textbf{20} (1998), 1049--1054.

\bibitem[Bar02a]{bartholdi:bsvlcs}
Laurent Bartholdi, \emph{The $2$-dimension series of the just-nonsolvable {BSV}
  group}, preprint, 2002.

\bibitem[Bar02b]{bartholdi:lpres}
Laurent Bartholdi, \emph{{$L$}-presentations and branch groups}, to appear in
  J. Algebra, 2002.

\bibitem[BGN02]{bartholdi-g-n:fractal}
Laurent Bartholdi, Rostislav~I. Grigorchuk, and Volodymyr~V. Nekrashevych,
  \emph{From fractal groups to fractal sets}, submitted, 2002.

\bibitem[BS97]{barnea-s:hausdorff}
Yiftach Barnea and Aner Shalev, \emph{{H}ausdorff dimension, pro-$p$ groups,
  and {K}ac-{M}oody algebras}, Trans. Amer. Math. Soc. \textbf{349} (1997),
  no.~12, 5073--5091.

\bibitem[BSV99]{brunner-s-v:nonsolvable}
Andrew~M. Brunner, Said~N. Sidki, and Ana~Cristina Vieira, \emph{A just
  nonsolvable torsion-free group defined on the binary tree}, J. Algebra
  \textbf{211} (1999), no.~1, 99--114.

\bibitem[CFP96]{cannon-f-p:thompson}
James~W. Cannon, William~J. Floyd, and William~R. Parry, \emph{Introductory
  notes on {R}ichard {T}hompson's groups}, Enseign. Math. (2) \textbf{42}
  (1996), no.~3-4, 215--256.

\bibitem[CGH99]{grigorchuk-h-s:paradox}
Tullio~G. {Ceccherini-Silberstein}, Rostislav~I. Grigorchuk, and Pierre~{de la}
  Harpe, \emph{Amenability and paradoxical decompositions for pseudogroups and
  discrete metric spaces}, Trudy Mat. Inst. Steklov. \textbf{224} (1999),
  no.~Algebra. Topol. Differ. Uravn. i ikh Prilozh., 68--111, Dedicated to
  Academician Lev Semenovich Pontryagin on the occasion of his 90th birthday
  (Russian).

\bibitem[Cho80]{chou:elementary}
Ching Chou, \emph{Elementary amenable groups}, Illinois J. Math. \textbf{24}
  (1980), no.~3, 396--407.

\bibitem[Day57]{day:amen}
Mahlon~M. Day, \emph{Amenable semigroups}, Illinois J. Math. \textbf{1} (1957),
  509--544.

\bibitem[Gri80]{grigorchuk:rw}
Rostislav~I. Grigorchuk, \emph{Symmetrical random walks on discrete groups},
  Multicomponent random systems, Dekker, New York, 1980, pp.~285--325.

\bibitem[Gri83]{grigorchuk:growth}
Rostislav~I. Grigorchuk, \emph{On the {M}ilnor problem of group growth}, Dokl.
  Akad. Nauk SSSR \textbf{271} (1983), no.~1, 30--33.

\bibitem[Gri98]{grigorchuk:amenEG}
Rostislav~I. Grigorchuk, \emph{An example of a finitely presented amenable
  group that does not belong to the class {E}{G}}, Mat. Sb. \textbf{189}
  (1998), no.~1, 79--100.

\bibitem[Gri00]{grigorchuk:jibg}
Rostislav~I. Grigorchuk, \emph{Just infinite branch groups}, New horizons in
  pro-$p$ groups (Markus P. F. du~Sautoy Dan~Segal and Aner Shalev, eds.),
  Birkh\"auser Boston, Boston, MA, 2000, pp.~121--179.

\bibitem[GS87]{ghys-s:thompson}
{\'E}tienne Ghys and Vlad Sergiescu, \emph{Sur un groupe remarquable de
  diff\'eomorphismes du cercle}, Comment. Math. Helv. \textbf{62} (1987),
  no.~2, 185--239.

\bibitem[GZ]{grigorchuk-z:torsionfree}
Rostislav~I. Grigorchuk and Anzdrzej \.Zuk, \emph{On a torsion-free weakly
  branch group defined by a three state automaton}, submitted.

\bibitem[HK97]{higson-k:bc}
Nigel Higson and Gennadi~G. Kasparov, \emph{Operator ${K}$-theory for groups
  which act properly and isometrically on {H}ilbert space}, Electron. Res.
  Announc. Amer. Math. Soc. \textbf{3} (1997), 131--142 (electronic).

\bibitem[Kes59]{kesten:rwalks}
Harry Kesten, \emph{Symmetric random walks on groups}, Trans. Amer. Math. Soc.
  \textbf{92} (1959), 336--354.

\bibitem[Mil68]{milnor:5603}
John~W. Milnor, \emph{Problem 5603}, Amer. Math. Monthly \textbf{75} (1968),
  685--686.

\bibitem[Ol{\cprime}80]{olshansky:invmean}
Alexander~Ju. Ol{\cprime}shanski{\u\i}, \emph{On the question of the existence
  of an invariant mean on a group}, Uspekhi Mat. Nauk \textbf{35} (1980),
  no.~4(214), 199--200.

\bibitem[OS01]{olshansky-s:nonamen}
Alexander~Yu. Ol{\cprime}shanski{\u\i} and Mark~V. Sapir, \emph{Non-amenable
  finitely presented torsion-by-cyclic groups}, Electron. Res. Announc. Amer.
  Math. Soc. \textbf{7} (2001), 63--71 (electronic).

\bibitem[R{\"o}v99]{rover:fpsg}
Claas~E. R{\"o}ver, \emph{Constructing finitely presented simple groups that
  contain {G}rigorchuk groups}, J. Algebra \textbf{220} (1999), no.~1,
  284--313.

\bibitem[Ser80]{serre:trees}
Jean-Pierre Serre, \emph{Trees}, Springer-Verlag, Berlin, 1980, Translated from
  the French by John Stillwell.

\bibitem[Sil01]{silva:phd}
Edm{\'e}ia Fernandes~da Silva, \emph{Uma fam\'\i lia de grupos quase
  n\~ao-sol\'uveis definida sobre \'arvores $n$-\'arias, $n\ge2$}, Ph.D.
  thesis, Universidade de Bras\'\i lia, 2001.

\bibitem[vN29]{vneumann:masses}
John von Neumann, \emph{Zur allgemeinen {Theorie} des {Masses}}, Fund. Math.
  \textbf{13} (1929), 73--116 and 333, = \emph{Collected works}, vol.\ I, pages
  599--643.

\end{thebibliography}
\end{document}